\documentclass[final,authoryear,11pt,singlecolumn]{elsarticle}
\usepackage[utf8]{inputenc}
\usepackage[english]{babel}
\usepackage{parskip}
\usepackage{setspace,longtable}
\usepackage{graphics}
\usepackage{enumerate}
\usepackage{amsmath, amssymb} 
\usepackage{verbatim}
\usepackage{xcolor}
\usepackage{url} 
\usepackage{array}
\usepackage[a4paper, inner=1.5cm, outer=1.5cm, top=2cm, bottom=3cm]{geometry}
\usepackage{algcompatible, algorithm}
\usepackage[utf8]{inputenc}
\usepackage[T1]{fontenc}
\usepackage{mathtools}
\usepackage[thinc]{esdiff}
\usepackage{soul}
\usepackage{hyperref, doi}
\usepackage{natbib} %reference in brackets

\setcitestyle{numbers}
\usepackage[title]{appendix}%for appendix title
\usepackage{ulem}
\usepackage{tabularx} %tablo sığdırma
\usepackage{makecell} %tablo sığdırma
\biboptions{round,sort}
\hypersetup{colorlinks = true, linkcolor = blue, citecolor = cyan}

\newcommand{\suchthat}{\;\ifnum\currentgrouptype=16 \middle\fi|\;}

\begin{document}
	
	\begin{frontmatter}
		
		\title{Course timetabling under conflict and lecturer preferences using bi-objective optimization model}
			
		\author[Yeditepe]{Murat Elhüseyni\corref{cor1}}
		\address[Yeditepe]{%
			Yeditepe University, Industrial Engineering, 34755, Istanbul, Turkey  }
		\ead{murat.elhuseyni@yeditepe.edu.tr}
				
		\cortext[cor1]{Corresponding author}
		
		\date{\today}

		\begin{abstract}

 In this work, we address a real-life course timetabling problem that incorporates several practical academic requirements, including the chronological precedence of theoretical sessions over practical ones, the enforcement of rest periods for lecturers, and the simultaneous assignment of lecturers and teaching assistants. The formulation involves a trade-off between minimizing course conflicts and maximizing schedule compactness by reducing the number of dispersed working days for lecturers under strict daily working hour constraints. To handle this trade-off, we develop a bi-objective mixed-integer linear programming (MILP) model solved using a weighted-sum objective approach. Computational experiments using a real-world dataset demonstrate that the proposed model reduces timetable conflicts by 13.33\%, improves lecturer schedule compactness by 26\%, and completely eliminates violations of theoretical–practical session precedence compared to the real-life timetable. Overall, the findings highlight the efficiency and robustness of the proposed MILP model for complex university timetabling problems.
		\end{abstract}
		
	\end{frontmatter}
	
	\begin{small}
		\begin{center}
		\end{center}
	\end{small}
	
	%	\maketitle
\section{Introduction}

Course timetabling is a well-known combinatorial optimization problem that arises in educational institutions when courses, instructors, students, rooms, and time slots must be coordinated under various constraints. The problem is challenging because numerous conflicting requirements must be satisfied simultaneously, including resource capacities, instructor availability, and student course combinations. As universities expand their programs and offer a wider variety of elective courses, the complexity of constructing feasible and high-quality timetables increases significantly.

To conduct university education seamlessly, departmental courses must be assigned appropriate times and locations. Hence, every university department constructs a course timetable to meet this requirement. The primary stakeholders of such timetables are students and lecturers, both of whom may be significantly affected by the resulting schedule. A closer examination of current practice reveals that timetables are generally created manually. The most common approach is to repeat the previous year’s timetable with minor adjustments. However, this practice may lead to several issues for stakeholders, such as additional course sections, night-owl students attending 8 a.m. classes, or professors teaching oversized sections because they are assigned to larger classrooms \citet{informs2022orms}. In addition, the ripple effect is highly prevalent in timetable construction; resolving a single scheduling conflict may require modifying many other course assignments. While visualization tools can help mitigate such difficulties \citet{informs2022orms}, they often require considerable time and effort, particularly when new elective courses are introduced or new lecturers join the department. 

Aside from the burden of the manual process, poorly designed timetables may create substantial inefficiencies for both students and instructors. The most important and delicate issue is \textit{conflict}. While decision makers generally prevent conflicts within the same cohort of students at the same grade level, conflicts frequently arise when a student wants to select a course from an upper or lower class in their program. This situation is not uncommon, as students may need to repeat a failed course or may take a higher-level course to reduce their workload in future semesters. Since course registration takes place after the departmental timetable is announced, this issue creates difficulties for both students and advisors who assist them during the registration period. In addition, lecturers typically prefer consolidated schedules rather than schedules dispersed across multiple days. The majority of courses consist only of theoretical lectures, while some courses also include practical sessions and therefore may have multiple sections. For each section, assigning practical sessions before the theoretical lectures can lead to pedagogical inconsistencies. Furthermore, even for lecturers who prefer compact schedules, long idle periods between the classes they teach may reduce their teaching efficiency. As a result, this study addresses the Course Timetable Scheduling problem drawn from a real-life setting, with the objectives of minimizing conflicts and improving instructors’ teaching efficiency while considering all of the conditions mentioned above. To facilitate and speed up the decision-making process, the problem is formulated as a \textit{Mixed Integer Linear Programming (MILP)} model in which institutional and operational constraints are explicitly incorporated.

To illustrate the structure of the timetable, we provide two tables, Table \ref{tab:sample} and Table \ref{tab:sampleinstructor}, which depict a sample timetable and instructor--assistant-based schedules for four time slots over a two-day period. The symbol $L$ represents practical sessions, whereas $T$ denotes theoretical sessions. If a course does not include a practical session, it is represented only by its course code. In this example, three courses are assigned to the first slot on Tuesday, resulting in two conflict violations. Overall, the timetable contains 10 violations. Furthermore, for the MATH207 course, the practical session precedes the theoretical session, which are scheduled on Monday and Tuesday, respectively, leading to a pedagogically improper assignment. Regarding Table \ref{tab:sampleinstructor}, lecturers are indicated after commas, whereas assistants assigned to practical sessions are shown in parentheses together with the $L$ symbol. The table also illustrates that instructor ME's courses are dispersed across two different days.

\begin{table}[H]
	\centering
	%\footnotesize
	\caption{A sample timetable}
	\label{tab:sample}
	\begin{tabular}{|c|c|c|}
		\hline
		\textbf{Period} & \textbf{Monday} & \textbf{Tuesday} \\ \hline
		
		1 & \makecell{MATH207(L1) \\ PHYS101(T)} & \makecell{MATH207(T) \\ TURK101 \\IE413} \\ \hline
		2 & \makecell{MATH207(L1) \\ PHYS101(T)} & \makecell{MATH207(T) \\ TURK101 \\ IE413} \\ \hline
		3 & \makecell{COME305 \\IE217} & \makecell{RCUL101 \\ COME312} \\ \hline
		4 & \makecell{COME305 \\IE217} & \makecell{ RCUL101 \\ COME312} \\ \hline
		
	\end{tabular}
\end{table}

\begin{table}[H]
	\centering
	%\footnotesize
	\caption{Department Instructor and Assistant based sample timetable }
	\label{tab:sampleinstructor}
	\begin{tabular}{|c|c|c|}
		\hline
		\textbf{Period} & \textbf{Monday} & \textbf{Tuesday} \\ \hline
		
		1 & MATH207 (L, SC) & \makecell{MATH207, ME \\ IE 413, BB} \\ \hline
		2 & \makecell{MATH207 (L, SC)} & \makecell{MATH207, ME \\ IE 413, BB} \\ \hline
		3 & \makecell{COME305, ME} & \makecell{RCUL101, HC \\ COME312, EB} \\ \hline
		4 & \makecell{COME305, ME} & \makecell{RCUL101, HC \\ COME312, EB} \\ \hline
		
	\end{tabular}
\end{table}

The remainder of this study is organized as follows. The next section presents the relevant literature. The methodology and model formulation are then described in detail, followed by computational results and analysis. The paper concludes with pedagogical insights and recommendations.

\section{Literature Review}

When examining the literature on the university timetabling problem, we categorize the studies according to their problem components and solution methodologies. We first discuss the instructor-related aspect of the problem. While some studies also consider instructor assignment to courses according to their capabilities and preferences \citep{arratia2021university, blanco1998course,boland2008new}, most of the literature assumes that course--instructor assignments are predetermined. Therefore, the distinction lies in how the studies handle instructor-related constraints and preferences. \citet{arratia2021university} address limitations on the maximum number of courses assigned to a professor and reserve time periods for administrative duties. \citet{navabi2022multicriteria} incorporate professors' health risk factors during the pandemic period. \citet{wang2003using}, \citet{ozturk2016kullanici}, \citet{yu2002genetic}, and \citet{palma2020considering} take teacher time-slot preferences into account. \citet{sorensen2014two} minimize instructors' idle times between lectures, whereas \citet{phillips2015integer} consider instructors' classroom preferences. \citet{daskalaki2004integer} handle both preferred teaching periods and preferences regarding practical sessions, and dispersion of teaching hours across different days of the week. Overall, only a limited number of studies address teaching assistant assignment, and to the best of our knowledge, no previous study explicitly considers consolidated instructor schedules and forces an idle period between consecutive courses as in our research.

The student aspect also plays an important role in university timetabling. \citet{yu2002genetic} and \citet{daskalaki2004integer} minimize idle periods between two active lectures. \citet{davison2025modelling} consider individual student preferences regarding course delivery modes and practical session-hour selection in a hybrid teaching environment. Similarly, \citet{barnhart2022course} and \citet{navabi2022multicriteria} address hybrid course mode selection together with individualized course registrations. \citet{boland2008new} propose a student-centered approach in which students select favorable conflict-free blocks to construct complete schedules. \citet{ceschia2012design} prohibit undesirable patterns such as attendance at the last event of the day, participation in more than two consecutive events, and having a single course on a given day. \citet{ozturk2016kullanici} and \citet{sorensen2014two} prevent students from having completely empty days. \citet{palma2020considering} address flexibility for delayed students by considering failure rates and course periodicity. Unlike this work, we indirectly address the issue by minimizing conflicts among courses from different grade levels. Moreover, our study enforces a pedagogical precedence relationship between theoretical and practical sessions.

The remaining aspects of the problem mainly involve multiple sections of opened courses \citep{arratia2021university,davison2025modelling,barnhart2022course,boland2008new,daskalaki2004integer,palma2020considering}, although this issue has been addressed only in a limited number of studies. In contrast, the majority of the literature focuses on classroom assignment \citep{davison2025modelling,barnhart2022course,boland2008new,ceschia2012design,daskalaki2004integer,deris2000timetable,navabi2022multicriteria,ozturk2016kullanici,yu2002genetic,zhang2005constructing,phillips2015integer,sorensen2014two,palma2020considering}. Although we do not address classroom assignment, our study is among the few works that simultaneously consider multiple sections and teaching assistant assignment alongside \citet{daskalaki2004integer}.

The university course timetabling problem has been widely studied using different modeling and solution paradigms. Early studies primarily relied on constraint satisfaction--based approaches, where feasibility was ensured by explicitly modeling institutional rules and resource limitations within constraint programming frameworks \citep{deris1999incorporating, deris2000timetable, blanco1998course, zhang2005constructing}. Due to the computational complexity of large-scale instances, researchers also developed metaheuristic approaches. For example, \citet{petrovic2004university} provide a survey of various metaheuristics applied to timetabling problems. \citet{yu2002genetic} and \citet{wang2003using} employ genetic algorithms, whereas \citet{mushi2010tabu} implement Tabu Search. Nevertheless, most studies formulate the problem as an integer programming model \citep{davison2025modelling,boland2008new,daskalaki2004integer,navabi2022multicriteria,ozturk2016kullanici,phillips2015integer}. Some studies decompose the problem into multiple stages and solve different subproblems sequentially \citep{blanco1998course, sorensen2014two}. Furthermore, many studies consider multiple objectives \citep{davison2025modelling,ceschia2012design,daskalaki2004integer,deris1999incorporating,navabi2022multicriteria,ozturk2016kullanici,wang2003using,phillips2015integer,sorensen2014two}. While most of them determine objective coefficients directly, only a few explicitly apply multi-objective optimization techniques. \citet{navabi2022multicriteria} employ hierarchical optimization, whereas \citet{phillips2015integer} and \citet{davison2025modelling} apply lexicographic optimization. In contrast, \citet{palma2020considering} transform multiple objectives into a single weighted-sum objective by normalizing them according to their best and worst values. Similar to these studies, our work also investigates the trade-off between objective function terms by varying their relative priorities.

Overall, the literature has not handled a study that simultaneously addresses assistant-supported practical sessions, multiple parallel sections, workload restrictions, enforcing an idle period between consecutive courses for an instructor, and compactness of the instructor schedule. In particular, there is no study that minimizes both conflict violations and fragmented instructor schedules while coordinating lecture and teaching assistant components within a bi-objective optimization framework. This gap motivates the present study, which models these aspects with the help of mixed-integer programming. 

\section{Our Approach and Contributions}

In this work, we derive a bi-objective MILP program to solve the course timetabling problem at a university department to facilitate decision-making. We consider many problem specific constraint while optimizing the tradeoff between number of conflicted courses and the number of dispersed courses for lecturers in the model. Finally, we compare the effectiveness of our proposed model against a timetable taken from real life. Our computational experiments suggest that the MILP model produces less conflicts, more consolidated schedule and better pedagogical assignment of theoretical and practical sessions over time. Our paper makes the following key contributions to the literature: 

\begin{itemize}
	\item We propose the first MILP model in the course timetabling literature that simultaneously minimizes course conflicts and lecturer course dispersion.  
	\item We investigate the tradeoff between these objectives using the weighted-sum method. 
	\item One of the few problems that handles lecturer and teaching assistant slot assignments at the same time. 
	\item We introduce the first course timetabling problem that enforces the precedence of theoretical sessions before practical ones while also forcing idle periods between consecutive courses of a lecturer.
	\item We introduce a new timetable dataset to the literature.
\end{itemize}

The rest of the paper is structured as follows: In Section~\ref{s:problemFormulation}, we give the details of the mathematical programming formulation. In Section~\ref{sec:comp}, we present the results of the computational experiment on the real life timetable dataset. Finally, in Section~\ref{s:conc}, we draw conclusions from  our paper with final remarks.

\section{Problem Formulation}
\label{s:problemFormulation}

In this case study, we consider the Fall 2021 course timetabling problem of the Industrial Engineering Department at Üsküdar University, where the author was previously affiliated. The timetable is prepared at the beginning of the academic term for a curriculum consisting of 9 periods distributed over 5 weekdays. The courses are divided into five categories: compulsory, departmental elective, field elective, social elective, and language elective. Students take the latter two categories from different departments. Faculty members teach most compulsory courses as well as departmental and field electives. However, some compulsory common courses are taught by lecturers from other departments. For instance, PHYS101 is taught by a lecturer from the Physics Department. Therefore, we do not consider timetable conflicts arising from courses taught by a lecturer outside the department. Moreover, some courses consist of multiple sections. MATH207 is an example, containing two sections, each divided into theoretical and practical sessions. From a pedagogical perspective, the theoretical session conducted by the lecturer should be followed by a practical session conducted by a teaching assistant. Hence, the problem requires assigning both lecturers and teaching assistants. In our modeling approach, we treat each section of MATH207 as a separate course. Since we do not possess sufficient practical session related data, such as capacities and shared usage by other departments, we exclude classroom assignment from the model. Table~\ref{tab:problemnotation} presents the sets, parameters, and decision variables used in the formulation.

\begin{table}[H]
	\caption{Sets, parameters, and decision variables}
%	\footnotesize
	\begin{center}
		\begin{tabular}{|p{4cm} p{10cm}|}
			\hline
			Sets & \\
			\hline
			$L$ & lecturers \\ 
			$A$ & assistants \\ 
			$C$ & courses \\
			$G$ & grades \\
			$K$ & types, $k=0$ theoretical, $k > 0$ practical \\
			$C_g$ & courses of grade $g$, $C_g \subset C$ \\
			$K_c$ & types of course $c$, $K_c \subset K$ \\
			$C_l$ & courses of lecturer $l$, $C_l \subset C$ \\ 
			$C_a$ & courses of assistant $a$, $C_a \subset C$ \\ 
			$D$ & days\\
			$T$ & time slots \\
			\hline            
			Parameters & \\
			\hline
			$D_{ck}$ & duration of type $k$ of course $c$ \\
			$\alpha$ & objective coefficient for the total number of conflicts \\
			\hline            
			Decision variables & \\ \hline
			$x_{ckdt} $ & 1, if course $c$ is assigned to type $k$ on day $d$ at time $t$, 0 otherwise\\
			$z_{ld}$ & 1, if lecturer $l$ teaches on day $d$, 0 otherwise \\
			$y_{ldt}$ & 1, if lecturer $l$ teaches on day $d$ at time $t$, 0 otherwise \\
			$as_{adt}$ & 1, if assistant $a$ teaches on day $d$ at time $t$, 0 otherwise \\
			$v_{dt}$ & the number of conflicts on day $d$ at time $t$ \\    
			$t_c$ & day of the theoretical course $c$ \\
			$p_{ck}$ & day of practical type $k$ of course $c$ \\  \hline
		\end{tabular}
	\end{center}
	\label{tab:problemnotation}
\end{table}

\vspace{-0.5cm}

%\subsection{Mathematical Programming Model}

Based on the notation defined in Table \ref{tab:problemnotation}, the bi-objective MILP model is formulated as follows:

\vspace{-1.5cm}

\begin{align}
	\min &\hspace{0.5em}   \alpha \sum_{d \in D} \sum_{t \in T} v_{dt}+ (1-\alpha) \sum_{l \in L} \sum_{d \in D} z_{ld}  \label{objfunc} \\     
	\mathrm{s.t.}   	
	&\hspace{0.5em} \sum_{i=0}^{i=3} y_{ld(t+i)}  \le 3 \quad & l\in L,d \in D, t= 1, \dots, |T|-3 \label{lectrest}    \\   
	%\textit{Lecturer-time assignment:}    
	&\hspace{0.5em} y_{ldt} = \sum_{c \in C_l} x_{c0dt}  \quad & l \in L,d \in D,t \in T \label{lectimeass} \\
	&\hspace{0.5em} as_{adt} = \sum_{c \in \mathcal{C}_a} x_{c1dt} \quad & a \in A, d \in D, t \in T \label{assisttimeass} 
\end{align}

\begin{align}
    	&\hspace{0.5em} \sum_{t \in T} y_{ldt} \leq 6 \quad & l \in L, d \in D \label{lecturerUB} \\
	&\hspace{0.5em} \sum_{t \in T} as_{adt} \leq 6  \quad& a \in A, d \in D \label{assistUB} \\
	&\hspace{0.5em} x_{c0dt} \le z_{ld} \quad & c \in C_l, l \in L,d \in D,t \in T \label{lecturedayLB} \\ %Lecturer-day relation 
    &\hspace{0.5em} z_{ld} \le \sum_{c \in C_l} \sum_{t \in T} x_{c0dt} \quad & l \in L,d \in D \label{lecturedayUB} \\
	&\hspace{0.5em}  d \cdot x_{c0dt} \le t_c \quad & c \in C,d \in D, t \in T \label{theprctc1} \\%courseass-practice time linkage 
	&\hspace{0.5em}  t_c \le d + (1 - x_{c0dt}) \cdot |D| \quad & c \in C,d \in D, t \in T \label{theprctc2} 
	\\
	&\hspace{0.5em}  d \cdot x_{ckdt} \le p_{ck} \quad & c \in C, k \in K_c, k>0, d \in D, t \in T \label{theprctc3} 
	\\ 
	&\hspace{0.5em}  p_{ck} \le d + (1 - x_{ckdt}) \cdot |D| \quad & c \in C, k \in K_c, k>0, d \in D, t \in T \label{theprctc4} 
	\\
	&\hspace{0.5em}  t_c + 1 \le p_{ck} \quad & c \in C, k \in K_c, k>0 \label{theprctconeday} %at least one day between them %note that it is for our teachers, outside is free
	\\       
	&\hspace{0.5em} \sum_{d \in D} \sum_{t \in T} x_{ckdt} = D_{ck} \quad & c \in C, k \in K_c \label{courseduration} %course duration must be satisfied 
	\\ 
	&\hspace{0.5em}  x_{ckdt} + x_{c'k'dt} \le 1 \quad &  c, c' \in C_g, c \neq c', k \in K_c, k' \in K_{c'}, d \in D,t \in T \label{conflictgradeprevent} %No conflict within same grade 
	\\  
	&\hspace{0.5em}  x_{ck55} = 0 \quad & c \in C, k \in K_c \label{fridayprayer} %friday prayer empty slot 
	\\ 
	&\hspace{0.5em}  \sum_{c \in C} \sum_{k \in K_c}  x_{ckdt} \le 1 + v_{dt} \quad & d \in D,  t \in T \label{conflictviol} %\textit{Violation definition:}
	\\ 
	&\hspace{0.5em}  x_{ckd0} \le x_{ckd1} \quad & c \in C, k \in K_c, D_{ck} \ge 2, d \in D \label{consecassgntzero} %\textit{Consecutive assignment t=0} 
	\\ 
	&\hspace{0.5em}  x_{ckdt} \le x_{ckd(t-1)} + x_{ckd(t+1)} \quad & c \in C, k \in K_c, D_{ck} \ge 2, d \in D, t=1,\dots,|T|-2 \label{consecassgntintrm} %\textit{Consecutive assignment interm} 
	\\ 
	&\hspace{0.5em}  x_{ckd(|T|-1)} \le x_{ckd(|T|-2)} \quad & c \in C, k \in K_c, D_{ck} \ge 2, d \in D \label{consecassgntlast} %\textit{Consecutive assignment last}
	\\
	&\hspace{0.5em} x_{ckdt} \in \{0,1\} & c \in C, k \in K_c, d \in D, t \in T \label{varbegin}\\
	&\hspace{0.5em} z_{ld} \in \{0,1\} & l \in L, d \in D \\
	&\hspace{0.5em} y_{ldt} \in \{0,1\} & l \in L, d \in D, t \in T \\
	&\hspace{0.5em} as_{adt}  \in \{0,1\} & a \in A, d \in D, t \in T \\
	&\hspace{0.5em} v_{dt} \ge 0 & d \in D, t \in T \\
	&\hspace{0.5em} t_c  \ge 0 & c \in C  \\
	&\hspace{0.5em} p_{ck} \ge 0 & c \in C, k \in K_c \label{varend}
\end{align}

%\begin{align}

%\end{align}

The objective function~\eqref{objfunc} minimizes the weighted sum of total conflict violations and the number of dispersed lecturer schedules. Constraint~\eqref{lectrest} ensures that lecturers have sufficient rest time between consecutive courses by prohibiting four consecutive hours of teaching. Constraints~\eqref{lectimeass}--\eqref{assisttimeass} link each theoretical course and practical session to a lecturer and an assistant, respectively. Constraints~\eqref{lecturerUB}--\eqref{assistUB} impose an upper limit of six working hours for lecturers and assistants, respectively. Constraints~\eqref{lecturedayLB}--\eqref{lecturedayUB} determine whether a lecturer teaches on a given day $d$. Constraints~\eqref{theprctc1}--\eqref{theprctc2} determine the day of the theoretical lecture, while constraints~\eqref{theprctc3}--\eqref{theprctc4} determine the day of the practical session. Constraints ~\eqref{theprctc1}-\eqref{theprctc4} link discrete assignment decisions to continuous variables and we call them \textit{pedagogical constraints.} Constraint~\eqref{theprctconeday} enforces that the theoretical session must precede the corresponding practical session. Constraint~\eqref{courseduration} ensures that course assignments exactly match the required course duration which is called \textit{completeness}. Constraint~\eqref{conflictgradeprevent} prevents conflicts among courses within the same grade. Constraint~\eqref{fridayprayer} enforces an empty time slot for Friday prayer. Constraint~\eqref{conflictviol} counts the number of conflicts occurring in each time period. Finally, constraints~\eqref{consecassgntzero}--\eqref{consecassgntlast} ensure that course sessions are assigned to consecutive time slots. Finally, constraints~\eqref{varbegin}-\eqref{varend} are variable domain restrictions.

\section{Computational Results}
\label{sec:comp}

In this section, we solve the problem in C++ using Visual Studio 2022 and CPLEX 22.1.1 solver. The experiments are conducted on a machine with two Intel(R) Xeon(R) Silver 4210R processors, 64 GB RAM, and 32 threads. We present the results of our mathematical programming for varying $\alpha$ in Table \ref{tab:performance_analysis}. 
As key performance indicators (KPI), $viol$ represents total conflict violation hours, $day$ denote the number days of dispersed lecturer schedules, and $time$ displays the CPU time. Once $\alpha=0$, the model only optimizes $day$ yielding (53,11) where the first term denotes $viol$ and the second one represents $day$. As is seen in the table, when $\alpha$ $\in$ [0.1,0.9], (52,11) occurred. When $\alpha=1$, optimizing only $viol$ results in (52,15). As a result, only one dominant solution emerges as (52,11). In terms of CPU time, $\alpha=1$ results in the least time whereas $\alpha=0.6$ results in the maximum time. The former one stems from focusing solely on one objective while the latter one arises from balancing prioritization between two terms. 

\begin{table}[H]
	\centering
	\caption{Key Performance Indicators with respect to $\alpha$ values}
	\label{tab:performance_analysis}
	\begin{tabular}{|c|c|c|c|}
		\hline
		\textbf{\textbf{$\alpha$}} & \textbf{viol} & \textbf{day} & \textbf{time} \\ \hline
		0   & 53 & 11 & 19 \\ \hline
		0.1 & 52 & 11 & 28 \\ \hline
		0.2 & 52 & 11 & 39 \\ \hline
		0.3 & 52 & 11 & 25 \\ \hline
		0.4 & 52 & 11 & 38 \\ \hline
		0.5 & 52 & 11 & 34 \\ \hline
		0.6 & 52 & 11 & 43 \\ \hline
		0.7 & 52 & 11 & 40 \\ \hline
		0.8 & 52 & 11 & 31 \\ \hline
		0.9 & 52 & 11 & 41 \\ \hline
		1   & 52 & 15 & 1  \\ \hline
	\end{tabular}
\end{table}

Among the (52,11) values, $\alpha=0.3$ yields the minimum computational time; hence, the detailed schedules are reported in Table \ref{tab:modeltimetable} and compared with the current department timetable in Table \ref{tab:currenttimetable}. The proposed model reduces the number of timetable conflicts from 60 to 52, corresponding to a 13.33\% improvement over the current schedule. Furthermore, there are 12 course-section pairs in the curriculum that require a theoretical session to precede a practical session. In the current real-life timetable, 7 of these 12 pairs violate this pedagogical constraint. The proposed model completely eliminates all 7 violations, achieving a 100\% compliance rate for pedagogical sequencing. 

At a finer level of detail, the maximum number of conflicts within a single time period decreases from 3 in Table 5 to 2 in Table 6. Regarding precedence violations in the current timetable, for instance, the theoretical and practical sessions of IE215 are scheduled on Tuesday in Section 1, while in Section 2 the practical session is assigned to Wednesday and the theoretical session to Thursday, thereby violating the required ordering. In Table 6, all practical sessions are scheduled after their corresponding theoretical sessions, thus fully satisfying the pedagogical constraints.

\begin{table}[H]
	\centering
	%\footnotesize
	\caption{Current Unified Timetable}
	\label{tab:currenttimetable}
	\begin{tabularx}{\textwidth}{|c|X|X|X|X|X|}
		\hline
		\textbf{Time} & \textbf{Monday} & \textbf{Tuesday} & \textbf{Wednesday} & \textbf{Thursday} & \textbf{Friday} \\ \hline
		
		1 & 
		\makecell{PHYS101(L1) \\ MATH203(T)} &
		MATH203(L) &
		RPSC109 &
		CHEM101(T) &
		ATA101 \\ \hline
		
		2 &
		\makecell{PHYS101(L1) \\ MATH203(T)} &
		MATH203(L) &
		\makecell{RPSC109 \\ MATH207(T1)} &
		\makecell{CHEM101(T) \\ IE215(T2)} &
		\makecell{ATA101 \\ MATH207(L2) \\ COME312 \\ IE413} \\ \hline
		
		3 &
		PHYS101(T) &
		\makecell{RCUL101 \\ IE215(T1)  \\ IE351(1)} &
		\makecell{RPSC109 \\ MATH207(T1)  \\ IE325(1)} &
		\makecell{CHEM101(L2) \\ IE215(T2) \\ SE305} &
		\makecell{TURK101 \\ MATH207(L2) \\ COME312 \\ IE413} \\ \hline
		
		4 &
		\makecell{PHYS101(T) \\ IE312} &
		\makecell{RCUL101 \\ IE215(T1) \\ IE351(1)} &
		IE325(1) &
		\makecell{CHEM101(L2) \\ IE211 \\ SE305} &
		TURK101 \\ \hline
		
		5 &
		\makecell{IE217(1) \\ IE312} &
		\makecell{IE215(L1) \\ OHS401}  &
		\makecell{MATH101(T) \\ MATH207(T2) \\ COME313} &
		\makecell{IE211 \\  IE305(T1) \\ IE323} &
		\\ \hline
		
		6 &
		IE217(1) &
		\makecell{ENG101 \\ IE215(L1) \\ IE301 \\ OHS401} &
		\makecell{MATH101(T) \\ MATH207(T2) \\ IE325(2) \\ COME313} &
		\makecell{CHEM101(L2) \\ IE223(T) \\ IE305(T1) \\ IE323} &
		\makecell{MATH207(L1) \\  IE305(L1) \\ IE413} \\ \hline
		
		7 &
		\makecell{CHEM101(L1) \\ IE217(2)} &
		\makecell{ENG101  \\ IE301} &
		IE325(2) &
		\makecell{CHEM101(L2) \\ IE223(T)} &
		\makecell{MATH207(L1) \\ IE305(L1) \\ IE413} \\ \hline
		
		8 &
		\makecell{CHEM101(L1) \\ IE217(2)} &
		\makecell{ENG101 \\ IE223(L1)  \\ IE351(2)} &
		\makecell{PHYS101(L2) \\ IE215(L2) \\ COME305} &
		\makecell{MATH101(L) \\ IE223(L2)  \\ IE305(T2)} &
		\makecell{IE223(L3) \\ IE305(L2)  \\ IE317} \\ \hline
		
		9 &
		&
		\makecell{IE223(L1) \\ IE351(2)} &
		\makecell{PHYS101(L2) \\ IE215(L2) \\ COME305} &
		\makecell{MATH101(L) \\  IE223(L2) \\ IE305(T2)} &
		\makecell{IE223(L3) \\  IE305(L2) \\ IE317} \\ \hline
		
	\end{tabularx}
\end{table}

\begin{table}[H]
	\centering
	%\footnotesize
	\caption{Unified Timetable of the model solution under $\alpha=0.3$}
	\label{tab:modeltimetable}
	\begin{tabularx}{\textwidth}{|c|X|X|X|X|X|}
		\hline
		\textbf{Time} & \textbf{Monday} & \textbf{Tuesday} & \textbf{Wednesday} & \textbf{Thursday} & \textbf{Friday} \\ \hline
		
		1 &
		\makecell{PHYS101(T) \\ COME305} &
		IE317 &
		IE217(2) &
		\makecell{IE305(T1) \\ IE413} &
		\makecell{CHEM101(L1) \\ MATH207(L1) \\ IE305(L1)} \\ \hline
		
		2 &
		\makecell{PHYS101(T) \\ COME305} &
		\makecell{PHYS101(L1) \\ IE215(L1) \\ IE317} &
		IE217(2) &
		\makecell{IE305(T1) \\ IE211 \\ IE413} &
		\makecell{CHEM101(L1) \\ MATH207(L1) \\ IE305(L1)} \\ \hline
		
		3 &
		IE223(T) &
		\makecell{PHYS101(L1) \\ IE215(L1)} &
		MATH207(T1)  &
		IE211 &
		IE223(L2) \\ \hline
		
		4 &
		\makecell{IE223(T) \\ SE305 \\ COME313} &
		\makecell{RCUL101 \\ IE215(L1) \\ IE325(1)} &
		\makecell{MATH207(T1) \\ MATH101(L) \\ IE351(1)} &
		\makecell{MATH207(L2) \\ COME312} &
		IE223(L2) \\ \hline
		
		5 &
		\makecell{CHEM101(T) \\ SE305 \\ COME313} &
		RCUL101 &
		\makecell{MATH101(L) \\ IE223(L1) \\ IE351(1)} &
		\makecell{MATH207(L2) \\ COME312} &
		\\ \hline
		
		6 &
		\makecell{IE217(1) \\ CHEM101(T)} &
		\makecell{MATH101(T) \\ IE215(T2) \\ IE325(2)} &
		\makecell{TURK101 \\ IE223(L1)} &
		\makecell{MATH203(L) \\ IE301} &
		\makecell{PHYS101(L2) \\ IE223(L3)} \\ \hline
		
		7 &
		\makecell{ENG101 \\ IE217(1) \\ IE323} &
		\makecell{IE215(T2) \\ IE325(2) \\ MATH101(T)} &
		\makecell{TURK101 \\ MATH207(T2)} &
		\makecell{RPSC109 \\ MATH203(L) \\ IE301} &
		\makecell{PHYS101(L2) \\ IE223(L3) \\ IE305(L2)} \\ \hline
		
		8 &
		\makecell{ENG101 \\ IE323 \\ OHS401} &
		\makecell{MATH203(T) \\ IE312} &
		\makecell{MATH207(T2) \\ IE351(2)} &
		\makecell{RPSC109 \\ IE305(T2)} &
		\makecell{IE305(L2) \\ CHEM101(L2) \\ IE215(L2)} \\ \hline
		
		9 &
		\makecell{ENG101 \\ OHS401} &
		\makecell{MATH203(T) \\ IE312} &
		\makecell{RPSC109 \\ IE351(2)} &
		\makecell{RPSC109 \\ IE305(T2)} &
		\makecell{CHEM101(L2) \\ IE215(L2)} \\ \hline
		
	\end{tabularx}
\end{table}

Concerning the instructor based timetables, we present the current and the model solution in Tables \ref{tab:instructortimetable} and \ref{tab:modeltimetable_instructor}, respectively. The results demonstrate that the number of working days is diminished from 15 to 11, so resulting in 26.67\% improvement. In detail, two people -OK and EB- work three days in Table \ref{tab:instructortimetable}. In Table \ref{tab:modeltimetable_instructor}, the number of maximum working days is reduced to 2, leading to more compact schedule. Besides, in Table \ref{tab:instructortimetable}, MA works 4 hours consecutively on Monday between 5 and 8 periods. This issue is resolved in Table \ref{tab:modeltimetable_instructor} in which no one works consecutively without resting. Regarding assistants, their practical sessions are assigned to the schedule without any conflict while satisfying the pedagogical precedence constraints.

\begin{table}[H]
	\centering
	%\footnotesize
    \small
	\caption{Current Department Instructor and Assistant based Timetable}
	\label{tab:instructortimetable}
	\begin{tabularx}{\textwidth}{|c|X|X|X|X|X|}
		\hline
		\textbf{Time} & \textbf{Monday} & \textbf{Tuesday} & \textbf{Wednesday} & \textbf{Thursday} & \textbf{Friday} \\ \hline
		
		1 &
		&
		&
		&
		&
		\\ \hline
		
		2 &
		&
		&
		MATH207(T1), ME &
		IE215, OK &
		\makecell{COME312, EB \\ IE413, BB \\ MATH207 (L2, SC)} \\ \hline
		
		3 &
		&
		\makecell{IE215, OK \\ IE351, MA} &
		\makecell{MATH207(T1), ME \\ IE325, EB} &
		\makecell{IE215, OK \\ SE305, ME} &
		\makecell{COME312, EB \\ IE413, BB \\ MATH207 (L2, SC)} \\ \hline
		
		4 &
		IE312, HC &
		\makecell{IE215, OK \\ IE351, MA} &
		IE325, EB &
		\makecell{IE211, MS \\ SE305, ME} &
		\\ \hline
		
		5 &
		\makecell{IE217(1), MA \\ IE312, HC} & IE215 (L1,SC)
		&
		\makecell{MATH207(T2), ME \\ COME313, OK} &
		\makecell{IE211, MS \\ IE305(T1), EB \\ IE323, OK} &
		\\ \hline
		
		6 &
		IE217(1), MA &
		\makecell{IE215 (L1,SC)  \\ IE301, BB} &
		\makecell{MATH207(T2), ME \\ IE325, EB \\ COME313, OK} &
		\makecell{ IE305(T1), EB \\ IE323, OK} &
		\makecell{IE413, BB \\ MATH207 (L1, SC) \\ IE305 (L1, TK)} \\ \hline
		
		7 &
		IE217(2), MA & IE301, BB &
		IE325, EB & &
		\makecell{IE413, BB \\ MATH207 (L1, SC) \\ IE305 (L1, TK)} \\ \hline
		
		8 &
		IE217(2), MA & IE351, MA &
		\makecell{IE215 (L2, SC) \\ COME305, ME} & IE305(T2), EB &
		\makecell{IE317, HC \\ IE305 (L2, TK)} \\ \hline
		
		9 &
		& IE351, MA &
		\makecell{IE215 (L2, SC) \\ COME305, ME} & IE305(T2), EB &
		\makecell{IE317, HC \\ IE305 (L2, TK)} \\ \hline
		
	\end{tabularx}
\end{table}

\begin{table}[H]
	\centering
	%\footnotesize
    \small
	\renewcommand\cellalign{tl}
	\caption{The model solution under $\alpha=0.3$ for Department Instructor and Assistant based Timetable}
	\label{tab:modeltimetable_instructor}
	
	\begin{tabularx}{\textwidth}{
			|c|
			>{\raggedright\arraybackslash}X|
			>{\raggedright\arraybackslash}X|
			>{\raggedright\arraybackslash}X|
			>{\raggedright\arraybackslash}X|
			>{\raggedright\arraybackslash}X|
		}
		\hline
		\textbf{Time} & \textbf{Monday} & \textbf{Tuesday} & \textbf{Wednesday} & \textbf{Thursday} & \textbf{Friday} \\ \hline
		
		1 &
		\makecell[l]{COME305, ME \\ IE215, OK} &
		IE317, HC &
		IE217(2), MA &
		\makecell[l]{IE305(T1), EB \\ IE413, BB} &
		\makecell[l]{IE305 (L1, TK) \\ MATH207 (L1, SC)} \\ \hline
		
		2 &
		\makecell[l]{COME305, ME \\ IE215, OK} &
		IE317, HC &
		IE217(2), MA &
		\makecell[l]{IE305(T2), EB \\ IE413, BB \\ IE211, MS} &
		\makecell[l]{IE305 (L1, TK) \\ MATH207 (L1, SC)} \\ \hline
		
		3 &
		&
		\makecell[l]{IE325, EB \\ IE215 (L1, SC)} &
		MATH207(T1), ME &
		IE211, MS &
		\\ \hline
		
		4 &
		\makecell[l]{SE305, ME \\ COME313, OK} &
		\makecell[l]{IE325, EB \\ RCUL101, HC \\ IE215 (L1, SC)} &
		\makecell[l]{MATH207(T1), ME \\ IE351, MA} &
		\makecell[l]{COME312, EB \\ MATH207 (L2, SC)} &
		\\ \hline
		
		5 &
		\makecell[l]{SE305, ME \\  IE217, MA} & RCUL101, HC &
		IE351, MA &
		\makecell[l]{COME312, EB  \\ MATH207 (L2, SC)} &
		\\ \hline
		
		6 &
		\makecell[l]{IE217(1), MA \\ COME305, ME} &
		\makecell[l]{IE215, OK \\ IE325, EB}&
		&
		IE301, BB & IE305 (L2, TK)  \\ \hline
		
		7 &
		\makecell[l]{IE217(1), MA \\ IE323, OK}&
		\makecell[l]{IE 325, EB \\ IE215, OK} &
		MATH207(T2), ME &
		IE301, BB & IE305 (L2, TK)  \\ \hline
		
		8 &
		IE323, OK & IE312, HC &
		\makecell[l]{MATH207(T2), ME \\ IE325, EB \\ IE217, MA \\ IE351, MA} & IE305(T2), EB  &
		\makecell[l]{IE215 (L2, SC)} \\ \hline
		
		9 &
		& IE312, HC &
		\makecell[l]{IE325, EB \\ IE217, MA \\ IE351, MA} & IE305(T2), EB  &
		\makecell[l]{IE215 (L2, SC)} \\ \hline
		
	\end{tabularx}
\end{table}

%\newpage
\section{Conclusion}
\label{s:conc}

In this study, we address the timetabling problem of Uskudar University for the Fall 2021 semester. Since the problem is bi-objective, we develop a weighted mixed-integer linear programming (MILP) model to obtain efficient schedules. Unlike previous studies in the literature, the proposed model simultaneously minimizes timetable conflicts and the number of dispersed working days of lecturers while enforcing pedagogical precedence constraints requiring practical sessions to be scheduled after their corresponding theoretical lectures. The behavior of the model is investigated under different objective coefficients, and the results associated with the best-performing coefficient are reported.

A detailed analysis of the results shows that all coefficient values smaller than 1 yield the same solution, indicating that the model is robust with respect to changes in the objective coefficients. Furthermore, the maximum computational time is observed when the coefficients are close to each other, which may be attributed to the solver attempting to balance the two competing objectives. Nevertheless, even the maximum solution time remains below 45 seconds, demonstrating that the problem can be solved to optimality efficiently.

Compared with the existing timetable, the proposed model provides substantial improvements. The reduction in timetable conflicts by more than 10\% enables a more convenient schedule for students, while nearly 30\% improvement in lecturer schedule compactness demonstrates clear superiority over the current timetable. The former improvement is mainly achieved by reducing the maximum number of conflicts within individual time slots, whereas the latter results from decreasing the maximum number of working days assigned to lecturers. In addition, although not explicitly included in the objective function, the model achieves a 100\% success rate in resolving all existing pedagogical precedence requirements together with minor improvements in resting times between consecutive courses. Finally, detailed instructor and assistant schedules are analyzed to further demonstrate the improvements achieved by the proposed approach.

As a continuation of this work, additional lecturer-oriented constraints may be incorporated. For example, the model could minimize idle gaps between consecutive courses or accommodate lecturer preferences regarding specific days off during the week. From the students’ perspective, idle time between consecutive courses on the same day could also be reduced to improve timetable compactness. Furthermore, schedules may be designed to provide students with a free day for individual study or personal activities. Finally, the proposed framework could be extended to better accommodate hybrid course delivery structures.

\section{Acknowledgements}

The author would like to thank Betul Sarı, Ceyda Macit, Seda Kurnaz, and Yusuf Serkan Yayıkçı for their efforts in preparing an earlier version of this manuscript.

\section*{Data Statement}

The codes are available at https://github.com/muratelhuseyni/Course-Timetabling

\bibliography{referencesapa}

\end{document}